\documentstyle{amsppt}
\magnification=\magstep1
\parindent=1em
\baselineskip 15pt
\hsize=12.3 cm
\vsize=18.5 cm
\NoRunningHeads
\pageno=1

\vsize=7.4in
\NoBlackBoxes


\def\tr{\text{tr }}

\def\diag{\text{diag }}



\def\Com{\text{Com }}

\topmatter

\title Spectral characterization of sums of commutators I\endtitle
\author
 N.J. Kalton
\endauthor
\address
Department of Mathematics,
University of Missouri,
Columbia, Mo.  65211, U.S.A.
\endaddress
\email  nigel\@math.missouri.edu
\endemail
\thanks The  author was supported  by  NSF
Grant DMS-9500125 \endthanks
\abstract

Suppose $\Cal J$ is a two-sided quasi-Banach ideal of compact operators
on
a separable infinite-dimensional Hilbert space $\Cal H$.  We show that an
operator
$T\in\Cal J$
can be expressed as finite linear combination of commutators $[A,B]$
where
$A\in\Cal J$ and $B\in\Cal B(\Cal H)$ if and only its eigenvalues
$(\lambda_n)$ (arranged in decreasing order of absolute value, repeated
according to algebraic multiplicity and augmented by zeros if necessary)
satisfy the condition that the diagonal operator
$\diag\{\frac1n(\lambda_1+\cdots +\lambda_n)\}$ is a member of
$\Cal J.$ This
answers (for quasi-Banach ideals) a question raised by Dykema, Figiel,
Weiss and Wodzicki.

\endabstract

\endtopmatter

\heading 1. Introduction \endheading

Let $\Cal H$ be a separable infinite-dimensional Hilbert space, and let
$\Cal J$ be a (two-sided) ideal contained in the ideal of compact
operators
$\Cal K(\Cal H)$ on $\Cal H.$
We define the {\it commutator subspace} $\Com \Cal J$ to be the closed
linear span of commutators $[A,B]=AB-BA$ where $A\in\Cal J$ and $B\in\Cal
B(\Cal H).$   It has been shown by Dykema, Figiel, Weiss and Wodzicki in
\cite{3} that if $\Cal I_1$ and $\Cal I_2$ are any two ideals then
the
linear span of commutators of the form $[A_1,A_2]$ where $A_j\in\Cal I_j$
for $j=1,2$ coincides with the commutator subspace $\Com\Cal J$ where
$\Cal J=\Cal I_1\Cal I_2.$

Pearcy and Topping (\cite {7}, cf. \cite{2}) showed that for the Schatten
ideal
$\Cal J=\Cal C_p$ when
$p>1$, we have $\Com\Cal C_p=\Cal C_p.$ They then raised the question
whether $\Com\Cal C_1=\{T\in\Cal C_1: \ \tr T=0\}.$  This question was
resolved negatively by Weiss \cite{8}, \cite{9}.  However, Anderson
\cite{1} showed
that in the case $p<1$ we have $\Com\Cal C_p=\{T\in\Cal C_p:\ \tr T=0\}.$

In \cite{6} a complete characterization of $\Com \Cal C_1$ was obtained.
It was shown that, in the case when $\Cal J=\Cal
C_1$
is the trace-class, then $T\in\Com\Cal C_1$ if and only if $T\in\Cal C_1$
and its eigenvalues $(\lambda_n(T))_{n=1}^{\infty}$, counted according to
algebraic multiplicity and  arranged in some order satisfying that
$(|\lambda_n(T)|)_{n=1}^{\infty}$ is decreasing, satisfies the inequality
$$
\sum_{n=1}^{\infty}\frac{|\lambda_1+\cdots+\lambda_n|}{n}<\infty.\tag
1.1$$
If the eigenvalue set of $T$ is finite one may extend the sequence
$\lambda_n(T)$ by including infinitely many zeroes.
This extended earlier partial results in \cite{8} and \cite{9}.

Since, for any ideal $\Cal J$, $\Com\Cal J$ is a self-adjoint subspace it
is clear that if
$T=H+iK$ is split into hermitian and skew-hermitian parts then
$T\in\Com\Cal J$ if and only if $H\in\Com\Cal J$ and $K\in\Com \Cal J.$
Thus to characterize $\Com\Cal J$ it is necessary only to characterize
the hermitian operators in $\Com\Cal J.$  In particular, the result above
shows that if $T\in\Cal C_1$ and the condition (1.1) is equivalent to the
pair of condition that $H$ and $K$ each satisfy (1.1).

Recently in \cite {3} a very general approach was developed which is
applicable to any ideal. It was shown that for any ideal
$\Cal
J$ an hermitian operator $H\in\Com\Cal J$ if and only if $H\in\Cal J$ and
the diagonal operator $\diag\{\frac1n(\lambda_1+\cdots+\lambda_n)\}$
belongs to $\Cal J$ where again $\lambda_n=\lambda_n(H)$ is the
eigenvalue sequence as above.  Although this yields an explicit test
for membership in $\Com\Cal J$ by the process of splitting into
hermitian and skew-hermitian parts,
 it leaves open the question whether same characterization in terms of
eigenvalues extends to all
operators in $\Com\Cal J,$ as in the case of the trace-class.

The aim of this paper is to show that for a fairly broad class of
``nice'' ideals the answer to this question is positive.  In a separate
note, in collaboration with Ken Dykema \cite{4}, we show that the
answer
is in general negative.

The condition we impose on an ideal $\Cal J$ is that it is {\it
geometrically stable}.  This means that if a diagonal operator
$\diag\{s_1,s_2,\ldots\}\in\Cal J$ where $s_1\ge s_2\ge \cdots\ge 0$ then
we also have $\diag\{u_1,u_2,\ldots\}\in\Cal J$ where $u_n=(s_1\ldots
s_n)^{1/n}.$  For any Banach or  quasi-Banach ideal (i.e. an ideal
equipped with an appropriate ideal quasi-norm) this condition is
automatic.

Let us note  that our results do depend in an essential way on the
results of \cite{3}, in that we use their result to reduce the
problem
to discussion of an operator of the form $T=H+iK$ where $H,K$ are
hermitian.

\demo{Acknowledgement}We are grateful to Ken Dykema
for
valuable discussions concerning the results of \cite{3} and this problem.
This work was done on a visit to Odense University and the author
would like to thank Niels Nielsen and the Department of Mathematics for
their hospitality.\enddemo

\heading 2. The key results\endheading

Let $T$ be a bounded operator on
a separable Hilbert space $\Cal H.$
 We denote  by  $s_n=s_n(T)$, for $n\ge 1$ the singular values of $T$.
It will be convenient to define $s_n$ for $n$ not an integer by
$s_n=s_{[n]+1}.$
With this notation we have the inequalities $s_n(S+T)\le
s_{n/2}(S)+s_{n/2}(T)$ and $s_n(ST)\le s_{n/2}(S)s_{n/2}(T).$

  If $T$ is compact we denote
 by $\lambda_n=\lambda_n(T)$ the eigenvalues of $T$
  repeated according to algebraic multiplicity and arranged in decreasing
  order                                                               of
   absolute value (this arrangement is not unique, so we
require some selection to be made).  Note that
$s_n(T)=|\lambda_n(|T|)|.$

Let us suppose that $f:\Bbb C\to \Bbb C$ is any function which vanishes
on a neighborhood of the origin.  Then we can define a functional
$\hat f:\Cal K(\Cal H)\to \Bbb C$ by the formula:
$$ \hat f(T) =\sum_{n=1}^{\infty}f(\lambda_n).$$

\proclaim{Lemma 2.1}(1) If $f$ is continuous then $\hat f$ is
continuous.\newline
(2) If $f$ is a Borel function then $\hat f$ is a Borel function.\newline
(3) If $f$ is continuous, real-valued and subharmonic then $\hat f$ is
plurisubharmonic on $\Cal K(\Cal H)$, i.e. if $S,T\in\Cal K(\Cal H)$ then
$$ \hat f(S) \le \int_0^{2\pi}\hat
f(S+e^{i\theta}T)\frac{d\theta}{2\pi}.$$\endproclaim

\demo{Proof}(1) Suppose $f(z)$ vanishes for $|z|\le \delta$ where
$\delta>0.$ If $T\in\Cal K(H)$ then suppose $m$ is the least integer
such that
$|\lambda_m(T)|<
\delta.$  Pick $\eta<\delta$ such that $|\lambda_m(T)|<\eta$ and if $m\ge
2$ then
$\eta<|\lambda_{m-1}(T)|$.  Now if $T_n$ is a sequence
of compact operators with $\lim_{n\to\infty}\|T_n-T\|=0$ then by results
in \cite{5} (see p.14, 18) we can find $n_0$ so that an ordering
$(\lambda'_k(T_n))_{k=1}^{\infty}$ of the eigenvalues of $T_n$ such that
$|\lambda'_k(T_n)|<\eta$ for $n\ge n_0$ and $k\ge m$ and
$\lim_{n\to\infty}\lambda'_k(T_n)=\lambda_k(T)$ if $k<m.$  It follows
easily that $\hat f(T_n)\to \hat f(T).$

(2) Observe that the set of $f$ such that $\hat f$ is Borel is closed
under pointwise convergence of sequences on $\Bbb C.$  (2) follows then
from (1).

(3) By (1) $\hat f$ is continuous and it therefore suffices to show
that
$$ \hat f(S) \le \int_0^{2\pi}\hat
f(S+e^{i\theta}T)\frac{d\theta}{2\pi},$$
for two finite rank operators $S,T.$  Hence we can suppose $S,T$ are
actually $n\times n$ matrices for some $n.$  Then the conclusion is
immediate from Proposition 5.2 of \cite{6}.\qed\enddemo

We now introduce certain functionals of the above type.
We define
$$\nu(T)=\sum_{|\lambda_n|\ge 1}1, \quad
\mu(T)=\sum_{n=1}^{\infty}\log_+ |\lambda_n|,\quad \text{ and }
 \chi(T)=\sum_{|\lambda_n|\ge 1}\lambda_n.$$
By the above lemma $\nu$ and $\chi$ are Borel functions while
$\mu$ is continuous.  If $|T|=(T^*T)^{1/2}$ the eigenvalues of $|T|$
correspond to the singular values of $T$. Notice that if $T$ is normal
then $\nu(T)=\nu(|T|).$

Notice also that each of the functionals $\mu,\nu$ and $\chi$ are
``disjointly additive'' in the sense that $\mu(S\oplus T)=\mu(S)+\mu(T)$
and etc.  Here $S\oplus T$ represents the operator defined on $\Cal
H\oplus \Cal H$ by $S\oplus T(x,y)=(Sx,Ty).$

\proclaim{Lemma 2.2}For any $T\in\Cal K(\Cal H)$ we have $0\le \mu(T)\le
\mu(|T|).$
\endproclaim

\demo{Proof}Suppose $\nu(T)=n.$  Then $\mu(T)=\log
|\lambda_1\ldots\lambda_n|^{1/n}\le \log |s_1\ldots s_n|^{1/n}\le
\mu(|T|)$  (see Gohberg-Krein p.37).\qed \enddemo

\proclaim{Lemma 2.3}If $S,T\in\Cal K(H)$ then $\nu(|S+T|) \le
\nu(2|S|)+\nu(2|T|).$  In particular if $T=H+iK$ with $H,K$ hermitian
then
$\nu(H) \le 2\nu(|T|)$.\endproclaim

\demo{Proof}These follow easily from the Weyl inequalities, that
$s_{m+n-1}(S+T)\le s_m(S)+s_{n}(T).$\qed\enddemo

\proclaim{Lemma 2.4}(1) If $T$ is a compact normal operator with
$T=H+iK$
for
$H,K$ hermitian then $|\chi(H)-\Re\chi(T)|\le \nu(T).$\newline
(2) If $T$ is any compact operator and $|\alpha|\le 1$ then
$|\alpha\chi(T)-\chi(\alpha T)| \le \nu(T).$\newline
(3) If $ (T_j)_{j=1}^n$ are compact
normal operators with
$T_1+\cdots+T_n=0$ then
$$ |\chi(T_1)+\cdots+\chi(T_n)| \le
(n-1)(\nu(T_1)+\cdots+\nu(T_n)).$$
\endproclaim

\demo{Proof}(1) We have $$\chi(H)-\Re\chi(T)=\sum_{\Re \lambda_n<1\le
|\lambda_n|}\Re\lambda_n.$$  The result follows immediately.

(2) We notice that
$$ \alpha\chi(T)-\chi(\alpha T) =\alpha \sum_{1\le
|\lambda_n|\le\alpha^{-1}}\lambda_n$$
whence the result follows.

(3) (Compare \cite {3}, Lemma 2.2.) Since each $T_j$ is normal there
exist self-adjoint projections
$P_j$ of rank $\nu(T_j)$ so that $\chi(T_j)=\tr (P_jTP_j)$ and
$\|P_j^{\perp}T_jP_j^{\perp}\|<1.$  Let $Q$ be a self-adjoint projection
of rank $d\le \sum_{j=1}^n\nu(T_j)$ whose range includes the range of
each $P_j.$  Then
$$ \tr (QT_jQ) =\tr (P_jT_jP_j) +\tr ((Q-P_j)T_j(Q-P_j))$$
and so
$$ \align
 |\sum_{j=1}^n\chi (T_j)|&= |\sum_{j=1}^n\tr(Q-P_j)T_j(Q-P_j)| \\
                     &\le \sum_{j=1}^n (d-\nu(T_j))\\
                     &\le (n-1)\sum_{j=1}^n\nu(T_j).\qed\endalign $$
\enddemo

Since $\chi$ is not a continuous function on $\Cal K(H)$ we will now
correct it to make a continuous function.  To this end we fix a
nondecreasing
$C^{\infty}-$function $\varphi:\Bbb R\to\Bbb R$ such that
$\varphi(x)=0$ if $x\le 0$, and $\varphi(x)=1$ if $x\ge 1.$
We define
$$ \chi_{\varphi}(T)=\sum_{n=1}^{\infty}\lambda_n\varphi(\log
\lambda_n)$$ where $\varphi(-\infty)=0.$   By Lemma 2.1, $\chi_{\phi}$ is
continuous.

\proclaim{Lemma 2.5}For any compact operator $T$, we have
$|\chi(T)-\chi_{\varphi}(T)|\le e\nu(T).$\endproclaim

\demo{Proof} $\chi(T)-\chi_{\varphi}(T)=\sum_{|\lambda_n(T)|\ge
1}(1-\varphi(\log|\lambda_n|)\lambda_n.$  Then result follows
immediately.\qed\enddemo

Now we define a second $C^2-$function $\psi:\Bbb R\to\Bbb R$ with the
properties that $\psi(x)=0$ if $x\le 0$ and
$\psi''(x)=e^x(|\varphi''(x)|+2|\varphi'(x)|)$ for $x\ge 0.$
$\psi$ is an increasing convex function, which is linear for $x\ge 1.$
Thus there is a constant $C_1>0$ such that $\psi(x)\le C_1\max(x,0)$ for
all
$x.$

We now proof a crucial lemma.

\proclaim{Lemma 2.6}Let $h:\Bbb C\to\Bbb R$ be defined by $h(0)=0$ and
$h(z)=\psi(\log|z|)-x\varphi(\log|z|) $ for $z\neq 0.$  Then $h$ is
subharmonic.
\endproclaim

\demo{Proof}Note that $h$ vanishes on a neighborhood of the origin. In
fact
$h$ is
$C^2$ so we check
$\nabla^2h.$ It is easy to check that (for $z\neq 0$),
$$\nabla^2(\psi(\log |z|))=|z|^{-2}\psi''(\log
|z|)=|z|^{-1}(|\varphi''(\log|z|)|+2|\varphi'(\log |z|)|).$$ We also have
 $$ \nabla^2(x\varphi(\log |z|))= 2\frac{x}{|z|^2}\varphi'(\log |z|)+
 \frac{x}{|z|^2}\varphi''(\log |z|).$$
Hence $\nabla^2h\ge 0.$\qed\enddemo

\proclaim{Theorem 2.7}Suppose $T\in\Cal K(\Cal H).$  Suppose $T=H+iK$
where $H=\frac12(T+T^*)$ and $K=\frac{1}{2i}(T-T^*).$   Then there is a
constant $C_2$ such that
$$ |\chi(H)- \Re \chi(T)|\le C_2\mu(2|T|)$$
and
$$ |\chi(K)-\Im \chi(T)|\le C_2\mu(2|T|).$$
\endproclaim

\demo{Proof}For convenience we will define the function
$F(z)=\frac12(T+zT^*).$ For $0\le\theta\le 2\pi,$ we have:
$$ F(e^{i\theta})-\frac{1+e^{i\theta}}2 H
-i\frac{1-e^{i\theta}
}{2}K=0.$$   Note that  each operator is normal and we also have from
Lemma 2.3 that
$\nu(|F(e^{i\theta})|)\le
2\nu(|T|)$ and $\nu(H),\nu(K)\le 2\nu(|T|).$
Appealing to Lemma 2.4 (3),
$$ |\chi(F(e^{i\theta}))-
\chi(\frac{1+e^{i\theta}}{2}H)
-\chi(i\frac{1-e^{i\theta}}{2}K)| \le 12\nu(|T|).$$

On the other hand Lemma 2.4 (2) gives that
$$ |\chi(\frac{1+e^{i\theta}}{2}H)-\frac{1+e^{i\theta}}{2}\chi(H)|\le
\nu(H)\le 2\nu(|T|)$$
and
$$|\chi(i\frac{1-e^{i\theta}}{2}K)-i\frac{1-e^{i\theta}}{2}\chi(K)| \le
2\nu(|T|).$$

Hence
$$ |\chi(F(e^{i\theta}))-
\frac{1+e^{i\theta}}{2}\chi(H)
-i\frac{1-e^{i\theta}}{2}\chi(K)| \le 16\nu(|T|).$$
Integrating over $\theta$ then gives
$$\left
|\int_0^{2\pi}\chi(F(e^{i\theta}))\frac{d\theta}{2\pi}-
\frac12(\chi(H)+i\chi(K))\right| \le 16\nu(|T|).$$

Taking real parts we have, in particular,
$$ \chi(H) \le 2\Re
\int_0^{2\pi}\chi((F(e^{i\theta}))\frac{d\theta}{2\pi}
+32\nu(|T|).$$

We now replace $\chi$ by the smoother function $\chi_{\varphi}$, and
using Lemma 2.5 with $e<3$:
$$ \chi(H) \le
2\int_0^{2\pi}\Re\chi_{\varphi}(F(e^{i\theta}))\frac{d\theta}{2\pi}
+44\nu(|T|).$$

Let $g(z)=\psi(\log |z|)$ for $z\neq 0$ and $g(0)=0.$
For any operator $S$ we can write
$$ \Re\chi_{\varphi}(S)=\hat g(S)-\hat h(S).$$  Note that $0\le\hat
g(S)\le C_1\mu(S).$  Thus
$$ \chi(H) \le
2C_1\int_0^{2\pi}
\mu(F(e^{i\theta}))\frac{d\theta}{2\pi}-
2\int_0^{2\pi}\hat
h(F(e^{i\theta}))\frac{d\theta}{2\pi}+44\nu(|T|).$$
Now since $h$ is subharmonic, the functional $\hat h$ is plurisubharmonic
by Lemma 2.1. Note that $\hat h(F(0))=\hat h(T/2)=\hat
g(T/2)-\Re\chi_{\varphi}(T/2). $   Hence, by 2.4(2) and 2.5,
$$
\align
2\int_0^{2\pi}\hat
h(F(e^{i\theta}))\frac{d\theta}{2\pi}&\ge 2\hat
g(T/2)-2\Re\chi_{\varphi}(T/2)\\
&\ge  -2\Re\chi(T/2)- 6\nu(T)\\
&\ge  -\Re\chi(T)-8\nu(T).\endalign
$$
Hence
$$ \chi(H) \le
2C_1\int_0^{2\pi}
\mu(F(e^{i\theta}))\frac{d\theta}{2\pi}+\Re\chi(T) +44\nu(|T|)+8\nu(T).$$
Note that for every $n$ we have
$s_n(F(e^{i\theta})) \le s_{n/2}(T)$ so that $\mu(F(e^{i\theta}))\le
2\mu(|T|).$
We thus can simplify, using Lemma 2.2, to
$$ \chi(H) \le \Re\chi(T) +4C_1\mu(|T|)+44\nu(|T|)+8\nu(T).$$
Now observe that
$\nu(T) \le (\log 2)^{-1}\mu(2T)$ so that for a suitable constant $C_2$
we have
$$ \chi(H)\le \Re\chi(T) + C_2\mu(2|T|).$$

We now consider $-T$, $iT$ and $-iT$ in place of $T$ and the theorem
follows.\qed\enddemo

\heading 3. The main results \endheading

Now suppose that $\Cal J$ is a two-sided ideal contained in $\Cal K(H)$.
We denote by $\Com \Cal J$ the linear subspace of $\Cal J$ generated by
all operators of the form $[S,T]=ST-TS$ for $T\in\Cal J$ and $T\in \Cal
B(\Cal H).$  It is shown in \cite{3} that if $\Cal I_1$ and $\Cal
I_2$ are ideals such that $\Cal I_1\Cal I_2=\Cal J$ then $\Com \Cal J$
coincides with the linear span of all $[S,T]$ where $S\in\Cal I_1$ and
$T\in\Cal I_2.$  It is clear that if $T=H+iK$ with $H,K$ hermitian then
$T\in\Com \Cal J$ if and only if $H,K\in\Com \Cal J.$
One of the main results of \cite{3} characterizes the hermitian
operators in
$\Com
\Cal J$.
We now state this result together with a useful rewording.

\proclaim{Theorem 3.1}Suppose $\Cal J$ is an ideal of compact operators
on $\Cal H.$  Let $N$ be a normal operator in $\Cal J,$ and let
$\lambda_n=\lambda_n(N).$ Then the following conditions on $N$ are
equivalent:\newline
(1) $N\in\Com \Cal J.$\newline
(2) $\diag\{\frac1n(\lambda_1+\cdots+\lambda_n)\}\in \Cal J.$ \newline
(3)  There exists $T\in\Cal J$ so that $\frac1n|\lambda_1+\cdots
+\lambda_n|\le s_n(T)$ for each $n\in \Bbb N.$\newline
(4)  There exists $T\in \Cal J$ such that for all $\alpha>0$ we have
$|\chi (\alpha N)| \le  \nu(\alpha |T|).$
\endproclaim

\demo{Proof}The equivalence of (1) and (2) for hermitian operators is
proved in
\cite{3}. We will first establish the equivalence of (2),(3) and (4).

(3) clearly implies (2).  We now check that (2) implies (3).
  For $m\ge n$ we
have $$\frac1m|\lambda_1+\cdots+\lambda_m|\le \max(\frac1n|\lambda_1+
\cdots +\lambda_n|, s_n(N)).$$
Let $T=\diag\{u_n\}$ where $u_n=\max_{m\ge
n}\frac1m|\lambda_1+\cdots+\lambda_m|.$  Then the above shows that
$T\in\Cal J$ and of course $\frac1n|\lambda_1+\cdots+\lambda_n|\le
s_n(T).$

Next we show (3) implies (4).   We can assume that
$T$ in
(3) satisfies
$s_n(T)\ge s_n(N)$ for every $n.$ Then if $\alpha^{-1}> s_1(T)$ we have
$\chi(\alpha N)=0$ and $\nu(T)=0.$ Otherwise let
$n$ be the largest integer such
$s_n(T)\ge \alpha^{-1}.$ Then $\nu(\alpha T)= n.$   Now suppose $m$ is
the largest integer so that $s_m(N)\ge \alpha^{-1}.$  Then $1\le m\le n$
and $\chi(\alpha N)=\alpha (\lambda_1+\cdots +\lambda_m).$

Thus we have
$$ \align |\chi(\alpha N)|&\le \alpha
|\lambda_1+\cdots+\lambda_{n+1}|+n+1\\
&\le (n+1)\alpha s_{n+1}(T) +n+1\\
& \le 2(n+1)\le 4\nu(\alpha T).\endalign
$$
This yields (4) with $T$ replaced by $T\oplus T\oplus T\oplus T.$

Now assume we have (4); we may assume $T$ is positive. We again assume
$s_n=s_n(T)\ge s_n(N)$ for all
$n.$  Now for any $n\in\Bbb N$ we have:
$$ \frac1n|\lambda_1+\cdots+\lambda_n| \le s_n+\frac1n|\sum_{|\lambda_k
|>s_n}\lambda_k|.$$
Suppose $\sigma>s_n$ is smaller than any $|\lambda_k|>s_n.$ Then
$$ \frac1n|\lambda_1+\cdots+\lambda_n| \le
s_n+\frac{\sigma}{n}|\chi(\sigma^{-1}N)|.$$
However $|\chi(\sigma^{-1}N)|\le \nu(\sigma^{-1}T)<n$, so that
$$ \frac1n|\lambda_1+\cdots+\lambda_n| \le
s_n+\sigma.$$
Letting $\sigma$ tend to $s_n$ yields
$$ \frac1n|\lambda_1+\cdots+\lambda_n| \le  2s_n$$
so that (3) holds if $T$ is replaced by $2T.$

Finally if $N=H+iK$ where $H,K$ are hermitian then we have by Lemma
2.4(1) that $|\chi(\alpha H)-\Re\chi(\alpha N)|,|\chi(\alpha K)-\Im
\chi(\alpha N)|\le \nu(\alpha T).$  Hence $N$ satisfies (4) if and only
if both $H$ and $K$ satisfy (4).  As remarked above, the results of
\cite{3} imply that for hermitian operators (1) and (2) and hence also
(1) and (4) are equivalent.  Thus (1) and (4) are also equivalent for
normal operators.
\qed\enddemo

Now let us introduce a stability condition on the ideal $\Cal J.$ We
will say that $\Cal J$ is {\it geometrically stable} if whenever
$\diag(s_1,s_2,\ldots)\in\Cal J$ with $s_1\ge s_2\ge\cdots$ then
$\diag (t_1,t_2,\ldots)\in\Cal J$ where $t_n=(s_1\ldots s_n)^{1/n}.$

We say that $\Cal J$ of compact operators is a quasi-Banach ideal (or
Schatten ideal) if it
can be equipped with a complete quasi-norm  $T\to \|T\|_{\Cal J}$ so that
we have the ideal property $\|ATB\|_{\Cal J}\le \|A\|_{\infty}\|T\|_{\Cal
J}\|B\|_{\infty}$ whenever $A,B\in\Cal B(\Cal H).$  Here we denote the
operator norm of $A$ by $\|A\|_{\infty}.$

\proclaim{Proposition 3.2}If $\Cal J$ is a quasi-Banach ideal then $\Cal
J$ is geometrically stable.\endproclaim

\demo{Proof}We can assume that for some $0<r\le 1$ that $\|\cdot\|_{\Cal
J}$ is an $r$-norm i.e. $\|S+T\|_{\Cal J}^r \le \|S\|_{\Cal
J}^r+\|T\|_{\Cal J}^r.$   Suppose $D=\diag (s_n)\in\Cal J,$ where $s_1\ge
s_2\ge\cdots.$  We recall our convention that $s_r=s_{[r]+1}$ where
$r> 0$ is not an integer. Then for each
$k\in\Bbb N$ we have $\|\diag (s_{n/2^k})\|_{\Cal J}\le
2^{k/r}\|D\|_{\Cal J}.$ Pick $\theta>k/r$.  Then by completeness the
series
$\sum_{k=0}^{\infty}2^{-\theta k}\diag (s_{n/2^k})$ converges in $\Cal
J.$  Thus $\diag (u_n)\in\Cal J$ where $u_n=\sum_{k=0}^{\infty}2^{-\theta
k}s_{n/2^k}.$ In fact $\|\diag (u_n)\|_{\Cal J}\le
(\sum_{k=0}^{\infty}2^{-(r\theta-k)})^{1/r}\|D\|_{\Cal J}.$
Now suppose $1\le j\le n.$  Pick $k\in \Bbb N$ so that $2^{-k}n\le j\le
2.2^{-k}n.$  Then $$s_j\le s_{n/2^k} \le 2^{k\theta}u_n\le
\left(\frac{2n}{j}\right)^{\theta}u_n.$$
Hence
$$ t_n \le 2^{\theta}n^{\theta}(n!)^{-\theta/ n} u_n\le Cu_n$$
for some constant $C.$  This implies $\diag(t_n)\in \Cal J$ and further
that $\|\diag(t_n)\|_{\Cal J}\le C\|D\|_{\Cal J}$ for some constant $C$
depending only on $r.$\qed\enddemo

We now prove the main result of this note, which, for the special case of
geometrically stable ideals, answers positively a question posed in
\cite{3}. It should be noted that in \cite{4} it is shown that for singly
generated ideals geometric stability is a necessary and sufficient
condition for the equivalence of (1) and (2) in Theorem 3.3.  Thus
in general (1) and (2) are not equivalent.

\proclaim{Theorem 3.3}Suppose $\Cal J$ is a geometrically stable  ideal
of compact operators
on $\Cal H$ (in particular this hold if $\Cal J$ is a quasi-Banach
ideal.) Let
$S\in\Cal J$ and let
$\lambda_n=\lambda_n(S).$ Then the following conditions on $S$ are
equivalent:\newline
(1) $S\in\Com \Cal J.$\newline
(2) $\diag\{\frac1n(\lambda_1+\cdots+\lambda_n)\}\in \Cal J.$ \newline
(3)  There exists $T\in\Cal J$ so that $\frac1n|\lambda_1+\cdots
+\lambda_n|\le s_n(T)$ for each $n\in \Bbb N.$\newline
(4)  There exists $T\in \Cal J$ such that for all $\alpha>0$ we have
$|\chi (\alpha S)| \le  \nu(\alpha |T|).$\newline
(5)  There exists $T\in\Cal J$ such that for all $\alpha>0$ we have
$|\chi (\alpha S)|\le \mu(\alpha |T|).$
\endproclaim

\demo{Proof}We first show that $N=\diag(\lambda_n)\in \Cal J.$  Indeed we
have $|\lambda_1\ldots\lambda_n|\le s_1\ldots s_n$ where $s_n=s_n(S).$
Thus $|\lambda_n|\le t_n=(s_1\ldots s_n)^{1/n}$ so that $N\in\Cal J$ by
geometric stability.  Hence (2),(3) and (4) are equivalent by Theorem
3.1.

It is clear that (4) implies (5)  (replace $T$ by $eT$).  Let us prove
that (5) implies (3).  We can suppose that $s_n=s_n(T)\ge s_n(S)$ for all
$n,$ and let $t_n=(s_1\ldots s_n)^{1/n}.$  Then
$$\align |\lambda_1+\cdots +\lambda_n| &\le |\sum_{|\lambda_k|\ge
s_n}\lambda_k| + ns_n\\
&\le s_n|\chi(s_n^{-1}S)|+ns_n\\
&\le s_n\mu(s_n^{-1}|T|) +ns_n\\
&\le s_n\sum_{k=1}^n\log(s_k/s_n)+ns_n\\
&\le ns_n\log(t_n/s_n) +ns_n.
\endalign
$$
   Hence  since $\log x\le x$
for all $x\ge 1,$
$$ \frac{|\lambda_1+\cdots+\lambda_n|}{n} \le t_n+s_n$$
and by the geometric stability of the ideal we have (3).

To conclude the proof we establish equivalence of (1) with (5).  To this
end note that if $S=H+iK$ with $H,K$ hermitian then by Theorem 2.7, there
is a constant $C_2$ so that, for $\alpha>0,$
$$|\Re\chi(\alpha S)- \chi(\alpha H)|,|\Im\chi(\alpha S)-\chi(\alpha
K)|\le
C_2\mu(2\alpha  |S|).$$

Now suppose first that $H,K$ both satisfy (5) so that there are operators
$T_1,T_2\in \Cal J$ with $|\chi(\alpha H)|\le \mu(\alpha T_1)$ and
$|\chi(\alpha K)|\le \mu(\alpha T_2)$ for $\alpha >0.$  Pick an integer
$n>2C_2$ and consider the operator $W=T_1\oplus T_2\oplus V$ where $V$ is
the direct sum of $n$ copies of $2|S|$.   Then $|\chi(\alpha S)|\le
\mu(\alpha W)$ for all $\alpha >0.$
Conversely if $S$ satisfies (5) for an appropriate operator $T$ then
$H$ and
$K$
satisfy (5) for $T$ replaced by $T\oplus V$
It follows $S$ satisfies (5) if and only if both $H$ and $K$ satisfy (5).
Now if $S\in\Com\Cal J$ then $H,K\in\Com\Cal J$ so that by Theorem 3.1
$H,K$ satisfy (2)-(4) and hence also (5).  Therefore (1) implies (5).

Conversely if (5) holds for S, then both $H,K$ satisfy (5) and hence also
(2)-(4); so by Theorem 3.1, $H,K\in\Com \Cal J$ and hence $S\in\Com\Cal
J$ i.e.
(5) implies (1).\qed\enddemo

\vskip10pt
\enddocument
\bye


\Refs


\ref\no{1}\by J.H. Anderson \paper Commutators in ideals of
trace-class
operators, II \jour Indiana Univ. Math. J. \vol 35\yr 1986\pages
373--378\endref

\ref\no{2}\by J.H. Anderson and L.N. Vaserstein
 \paper Commutators in ideals of trace-class
operators \jour Indiana Univ. Math. J. \vol 35\yr 1986\pages
345-372\endref

\ref\no{3}\by K.J. Dykema, T. Figiel, G. Weiss and M. Wodzicki \paper
The commutator structure of
 operator ideals \paperinfo preprint \yr 1997 \endref

\ref\no{4}\by K.J. Dykema and N.J. Kalton \paper Spectral
characterization of
sums of commutators II\paperinfo preprint \yr 1997\endref

\ref\no{5}\by I.C. Gohberg and M.G. Krein \book Introduction to the
theory
of linear nonselfadjoint operators \bookinfo Translations of Mathematical
Monographs \vol 18 \publ Amer. Math. Soc. \publaddr Providence \yr 1969
\endref

\ref\no{6}\by N.J. Kalton \paper Trace-class operators and commutators
\jour J. Functional Analysis \vol 86 \yr 1989 \pages 41--74\endref


\ref\no{7}\by C.M. Pearcy and D. Topping \paper On commutators of ideals
of compact operators \jour Michigan Math. J. \vol 18\yr 1971 \pages
247--252\endref

\ref\no{8}\by G. Weiss\paper Commutators of Hilbert-Schmidt operators
II
\jour
Integral Equns. Operator Theory \vol 3 \yr 1980 \pages 574-600 \endref

\ref\no{9}\by G. Weiss\paper Commutators of Hilbert-Schmidt operators
I
\jour Integral Equns. Operator Theory \vol 9 \yr 1986 \pages 877-892
\endref



\endRefs
\bye